\input amstex
\documentstyle{amsppt}

\magnification=1200

\NoBlackBoxes

\def\R{\Bbb R}
\def\P{\Cal P}
\def\Q{\Cal Q}
\def\c{\colon}
\def\la{\langle\,}
\def\ra{\,\rangle}
\def\Vp{V^*\oplus V}
\def\tht{\thetag}
\def\gln{\frak{gl}(n)}
\def\Ugln{\Cal U(\gln)}
\def\Zgln{\frak Z(\gln)}
\def\M{\operatorname{M}}
\def\id{\operatorname{id}}
\def\sgn{\operatorname{sgn}}
\def\tr{\operatorname{tr}}
\def\Tab{\operatorname{Tab}}
\def\a{\alpha}
\def\p{\partial}
\def\si{\sigma}
\def\Si{\Sigma}
\def\d{\delta}
\def\l{\lambda}
\def\nea{\nearrow}
\def\RSk{\R[S(k)]}
\def\pst{\psi_{TT'}}

\def\Pst{\Psi_{TT'}}

\def\ov{\overset}
\def\G{\Gamma}
\def\chm{\chi^\mu}
\def\ovr{\,/\,}
\def\ET{{\Bbb E}_T}
\def\ETT{{\Bbb E}_{T T'}}

\def\A{\Bbb A}
\def\Sm{{\Bbb S}_{\mu}}
\def\gr{\operatorname{gr}}
\def\sm{s^*_{\mu}}
\def\f{\downharpoonright}

\def\Ls{\Lambda^*}
\def\gi{(g^{-1})}

\def\Ad{\operatorname{Ad}}

\def\pairfill{{$\mathsurround=0pt
\bracelu\leaders\vrule\hfill
\braceru$}}

\def\lpairfill{{$\mathsurround=0pt
\bracelu\leaders\vrule\hfill
{}$}}

\def\rpairfill{{$\mathsurround=0pt
{}\leaders\vrule\hfill
\braceru$}}

\def\pair#1{
\mathop{\vtop{\ialign{##\crcr
$\hfil\displaystyle{#1}\hfil$\crcr
\noalign{\kern3pt\nointerlineskip}
\kern3pt\pairfill\kern3pt
\crcr\noalign{\kern3pt}}}}\limits}

\def\chainpair#1{
\mathop{\vtop{\ialign{##\crcr
$\hfil\displaystyle{#1}\hfil$\crcr
\noalign{\kern3pt\nointerlineskip}
\kern-5pt\pairfill\kern3pt
\crcr\noalign{\kern3pt}}}}\limits}

\def\lchainpair#1{
\mathop{\vtop{\ialign{##\crcr
$\hfil\displaystyle{#1}\hfil$\crcr
\noalign{\kern3pt\nointerlineskip}
\kern-5pt\lpairfill\kern3pt
\crcr\noalign{\kern3pt}}}}\limits}

\def\rpair#1{
\mathop{\vtop{\ialign{##\crcr
$\hfil\displaystyle{#1}\hfil$\crcr
\noalign{\kern3pt\nointerlineskip}
\kern3pt\rpairfill\kern3pt
\crcr\noalign{\kern3pt}}}}\limits}

\def\dpair#1{
\mathop{\vtop{\ialign{##\crcr
$\hfil\displaystyle{#1}\hfil$\crcr
\noalign{\kern-1pt\nointerlineskip}
\kern3pt\pairfill\kern3pt
\crcr\noalign{\kern3pt}}}}\limits}

\def\dchainpair#1{
\mathop{\vtop{\ialign{##\crcr
$\hfil\displaystyle{#1}\hfil$\crcr
\noalign{\kern-1pt\nointerlineskip}
\kern-5pt\pairfill\kern3pt
\crcr\noalign{\kern3pt}}}}\limits}

\topmatter
\title
Young basis, Wick formula, \\
and higher Capelli identities 
\endtitle
\author
Andrei Okounkov
\endauthor
\affil
Institute for Problems of Information Transmission,
Moscow, and\\
Institute for Advanced Study, Princeton
\endaffil
\abstract
We prove Capelli type identities which
involve the whole universal 
enveloping algebra $\Ugln$ and matrix elements
of irreducible representations of the symmetric group. 
These identities generalize higher
Capelli identities for the center of $\Ugln$
introduced in the author's paper \cite{Ok}.
The main role in the proof play the Jucis-Murphy
elements.
\endabstract 
\address
School of Mathematics, Institute for Advanced Study,
Princeton, NJ 08540
\endaddress
\email
okounkov\@math.ias.edu
\endemail
\thanks
The author was supported by International Science Foundation 
and by Russian Basic Research Foundation. During his stay
at the Insitute for Advanced Study the author was supported
by the NSF grant DMS 9304580.
\endthanks
\endtopmatter

\head 
1. Introduction
\endhead

\subhead 1.1\endsubhead
Identify the  the standard generators $E_{ij}$ of the
Lie algebra with the following vector fields on
the vector space $\M(n,m)$ of $n\times m$ matrices
$$
E_{ij}=\sum_\alpha x_{i\alpha} \p_{j\alpha} \,, \tag 1.1
$$
where $x_{ij}$ are the natural coordinates in $\M(n,m)$ and
$\p_{ij}$ are the dual partial derivatives.
Suppose $m\ge n$. Identify the universal enveloping algebra $\Ugln$
with the algebra of differential operators
with polynomial coefficients on $\M(n,m)$ 
invariant under the action of $GL(m)$ by
right multiplication.

We obtain an explicit expression for a very large
and remarkable family of right-invariant differential
operators on $\M(n,m)$ in terms of generators $E_{ij}$.
In particular, our result is a generalization of
the following classical Capelli identity \cite{C}
$$
\multline
\sum_{1\le i_1<i_2<\dots<i_k\le n}
\operatorname{row-det}
\left(
\matrix
\format \l\quad &\l\quad &\l\quad &\l \\
E_{i_1i_1} & E_{i_1i_2} & \hdots & E_{i_1i_k} \\
E_{i_2i_1} & E_{i_2i_2}+1 &  & E_{i_2i_k} \\
\vdots & & \ddots & \vdots \\
E_{i_ki_1} & E_{i_ki_2} & \hdots & E_{i_ki_k}+k-1 
\endmatrix
\right)
= \\
\sum_{1\le i_1<i_2<\dots<i_k\le n} \quad
\sum_{1\le j_1<j_2<\dots<j_k\le m}
\det(x_{i_aj_b})_{1\le a,b,\le k} \,
\det(\p_{j_ai_b})_{1\le a,b,\le k} \,, 
\endmultline \tag 1.2
$$
where $k=1,2,\dots,n$ and 
the row-determinant of a matrix $A=(A_{ij})$ with
entries in a non-commutative algebra $\A$ is defined
by the following formula
$$
\operatorname{row-det} A =
\sum_{s\in S(n)} \sgn(s)\, A_{1\,s(1)} 
A_{2\,s(2)} \dots A_{n\,s(n)}\,.
$$
A detailed discussion of this famous result of
the classical invariant theory can be found in \cite{HU,KS}.

\subhead 1.2\endsubhead
Introduce some notation. It is convenient
to use matrices with entries in a non-commutative
algebra. Let $E$, $X$, and $D$ denote
matrices with entries $E_{ij}$, $x_{ij}$ and $\p_{ij}$
respectively. Then \tht{1.1} is equivalent to
$$
E=X\cdot D'\,,
$$
where prime stands for transposition. Introduce also the 
matrix
$$
E-u=(E_{ij}-u\cdot\d_{ij})_{ij}\,,
$$
which depends on a formal parameter $u$.

A $n\times n$ matrix $A$
with entries $A_{ij}$ in a non-commutative algebra $\A$ can
be considered as an element
$$
A=\sum_{ij} A_{ij} \otimes e_{ij}  \,\in\, \A\otimes \M(n) \,,
$$
where $e_{ij}$ are standard matrix units in $\M(n)$. The
tensor product of two such matrices $A$ and $B$ is defined by
$$
A\otimes B = \sum_{i,j,k,l} A_{ij} B_{kl}
\otimes e_{ij} \otimes e_{kl} \, \in  \A\otimes \M(n)^{\otimes 2} \,.
$$
Define the trace of an element
of $\A\otimes \M(n)^{\otimes n}$ by
$$
\tr \big( \sum_
{i_1,j_1,\dots,i_n,j_n}
A_{i_1,j_1,\dots,i_n,j_n} \otimes
e_{i_1,j_1}\otimes\dots\otimes e_{i_n,j_n} \big) =
\sum_{i_1\dots,i_n}
A_{i_1,i_1,\dots,i_n,i_n} \, \in \A \,.
$$
The symmetric group $S(k)$ acts in the vector space of $k$-tensors, 
so that we have a representation
$$
S(k)\to\M(n)^{\otimes k}\,.
$$

Let $\mu$ be a Young diagram with $k$ boxes. Let $V^\mu$ be the
corresponding irreducible $S(k)$-module and let $\chm$ be its
character. Consider $\chm$ as an element of the group algebra
of $S(k)$
$$
\chm = \sum_{s\in S(k)} \chm(s) \cdot s \, \in \RSk \,.
$$
Let $T$ be a Young tableau of shape $\mu$ and let $v_T$
be the corresponding vector in the Young basis for $V^\mu$.
(We recall some basic facts from  representation theory
of $S(k)$ in section 3 below.) Let $T$ and $T'$ be two
Young tableaux of shape $\mu$. Consider the following matrix
element
$$
\Pst=\sum_{s\in S(k)}
(s\cdot v_T , v_{T'})\cdot s^{-1} \, \in \RSk \,.
$$
Let $\a=(i,j)$ be a box from $\mu$. The number
$$
c(\a)=j-i
$$
is called the content of the box $\a$. Write $c_T(i)$
for the content of the box number $i$ in 
the tableau $T$.

\subhead 1.3\endsubhead
It is not difficult to see (see \cite{MNO} or below) that the
Capelli identity can be restated as follows
$$
\tr \left(
E\otimes (E+1)\otimes \dots \otimes (E+k-1)
\cdot \chi^{(1^k)} \right) = \tr \left( X^{\otimes k} \cdot
(D')^{\otimes k} \cdot \chi^{(1^k)} \right)\,.
$$
The main result we prove in this paper is the following
identity
\proclaim{Main Theorem} Let $T$ and $T'$ be two Young tableaux
of the same shape. Then
$$
(E-c_T(1))\otimes \dots \otimes (E-c_T(k))
\cdot \Pst =  X^{\otimes k} \cdot
(D')^{\otimes k} \cdot \Pst \,. \tag 1.3
$$
\endproclaim

The proof is based of some remarkable properties of
Jucys-Murphy elements \tht{3.2}.

This matrix identity is equivalent to the following  $n^{2k}$
identities: for any two $k$-tuples of indexes
$$
i_1,\dots,i_k,\quad j_1,\dots,j_k
$$
we have the equality of the corresponding matrix elements
$$
\multline
\sum_{s\in S(k)} (s\cdot v_T,v_{T'})
(E_{i_1j_{s(1)}}-c_T(1) \d_{i_1j_{s(1)}})  \dots
(E_{i_kj_{s(k)}}-c_T(k) \d_{i_kj_{s(k)}})  = \\
\sum_{s\in S(k)} (s\cdot v_T,v_{T'})
\sum_{\a_1,\dots,\a_k}
x_{i_1\a_1} \dots x_{i_k\a_k} \,
\p_{j_{s(1)}\a_1} \dots \p_{j_{s(k)}\a_k} \,.
\endmultline \tag 1.4
$$
In contrast to Capelli identity, these identities involve
not only the generators of the center $\Zgln$ of $\Ugln$.
It is easy to see that linear combinations of \tht{1.3}
span the whole algebra $\Ugln$. Moreover, it is easy to
see that there are much more identities \tht{1.3} than
linearly independent elements of $\Ugln$. 

Observe also that the identity \tht{1.3} is linear in $v_{T'}$;
therefore this vector can be replaced by an arbitrary
vector in $V^\mu$.

The  matrix
$$
\ETT=(E-c_T(1))\otimes \dots \otimes (E-c_T(k))
\cdot \Pst 
$$
should be called, perhaps, the {\it fusion} of $k$ matrices $E$
as its structure is the same as the fusion of $R$-matrices,
see \cite{Ch,KuSR,KuR}. The identity \tht{1.3} is one of the
interesting properties of this matrix.

In section 5 we specialize the identities for central
elements of $\Ugln$ by taking trace. We recover
higher Capelli identites introduced in the author's
paper \cite{Ok}. We have to mention that the proof of \tht{1.3}
given below is direct and does not require
$R$-matrix formalism used in \cite{Ok}.

The element
$$
\Sm=(\dim\mu\ovr k!\,) \tr{\Bbb E}_{TT}
$$
depends only on $\mu$, not on $T$; it was called in \cite{Ok} the
{\it quantum $\mu$-immanant}. Quantum immanants form
a very destinguished linear basis in $\Zgln$. We
recall only some basic facts about them from
\cite{Ok} and \cite{OO}. 

I am grateful to S.~Khoroshkin, V.~Ginzburg,
S.~Kerov, A.~Vershik, M.~Semenov-Tian-Shansky and
especially to  M. Nazarov and G. Olshanski for
helpful discussions. M. Nazarov proved a statement
equivalent to our main theorem using R-matrix methods.
His proof will be published in \cite{N}.


\head
2. Wick formula.
\endhead

\subhead 2.1 \endsubhead
Let $V$ be a vector space. By $D(V)$ denote the
algebra of differential operators on $V$ with
polynomial coefficients. The algebra $D(V)$ is
generated by constant vector fields $v\in V$ and
by multiplications by linear functions $\xi\in V^*$
subject to Heisenberg commutation relations
$$
[v,\xi]=\la\xi,v\ra\,,
$$
where $\la,\ra$ is the canonical pairing
$$
V^*\otimes V\to\R\,.
$$

By $S(\Vp)$ denote the symmetric algebra of the
vector space $\Vp$. Introduce the following linear
isomorphism
$$
S(\Vp)@>\quad\dsize{::}\quad>>D(V)\,,
$$
called {\it normal ordering}.
By definition, the normal ordering places all
multiplications by functions to the left and all
constant vector fields to the right; for example,
$$
\c\xi_1v_1v_2\xi_2\xi_3v_3\c\,=
\xi_1\xi_2\xi_3v_1v_2v_3\,,
\qquad \xi_i\in V^*,v_i\in V\,.
$$
By definition, put
$$
\pair{A\,B}=AB-\c AB\c,\qquad A,B\in\Vp\,. \tag 2.1
$$
This is a number; it is linear in $A$ and $B$. Clearly,
$$
\align
\pair{v\,\xi}&=\la\xi,v\ra\,,\\
\pair{\xi\,v}&=0\,,\qquad\xi\in V^*,v\in V\,.
\endalign
$$
Let the pairing
$$
\dots\pair{A\dots B}\dots
$$
mean that this pair should be replaced by the number \tht{2.1}.
The following theorem can be easily proved by induction.

\proclaim{\smc Theorem (Wick)}
Suppose $A_1,\dots,A_k\in\Vp\subset D(V)$. Then
$$
\align
A_1\dots A_k=&
\c A_1\dots A_k\c +\\
&\sum_{1\le i<j\le k} 
\c A_1\dots\pair{A_i\dots A_j}\dots A_k\c +\\
&\sum_{i,j,p,q}
\c A_1\dots\pair{A_i\dots A_j}
\dots\pair{A_p\dots A_q}\dots A_k\c + \dots \,,
\endalign
$$
where the sum is over all possible pairings in the set
$\{1,\dots,k\}$.
\endproclaim
For example,
$$
A_1A_2A_3=\c A_1A_2A_3\c + \c\pair{A_1A_2}A_3\c+
\c A_1\pair{A_2A_3}\c+\c\pair{A_1A_2A_3}\c\,.
$$

\subhead 2.2 \endsubhead
Recall that we identify $\Ugln$ with the algebra of right-invariant
differential operators on $\M(n,m)$ with polynomial coefficients.

Consider the following linear isomorphism
$$
S(\gln)@>\,\si\,>>\Ugln\,
$$
introduced by G.~Olshanski in \cite{Ol1}; in \cite{KO} it was
called the {\it special} symmetrization. The definition of
$\si$ is equivalent to the following (see lemma 2.2.12 in \cite{Ol1})
$$
\si(E_{i_1j_1}\dots E_{i_kj_k})
=\sum_{\a_1,\dots,\a_k}
x_{i_1\a_1}\dots x_{i_k\a_k}
\p_{j_1\a_1}\dots \p_{j_k\a_k}\,. \tag 2.2
$$
It is easy to see that the RHS of \tht{2.2} is a right-invariant
differential operator and hence an element of $\Ugln$. By
analogy to the normal ordering let us call the map $\si$
the {\it normal} symmetrization and denote it by colons
$$
\c E_{i_1j_1}\dots E_{i_kj_k}\c =
\si(E_{i_1j_1}\dots E_{i_kj_k}) \,.
$$

Suppose $A,B\in\gln$. Put
$$
\pair{A\,B}=AB-\c AB\c \quad \in\gln\,.
$$
It is easy to see that this is simply the matrix multiplication
$$
\pair{E_{ij}E_{pq}}=\d_{jp}E_{iq}\,.
$$
Observe that chain pairings like
$$
\c A_1\dots\pair{A_a\dots A_b}\chainpair{\dots A_c} 
\dots A_k\c\,,\tag 2.3
$$
where the end of a brace is at the same time the
beginning of a new brace, make perfect sence in the 
case of $\gln$. The pairing \tht{2.3} simply means
that the three matrices should be replaced by their
matrix product. The following theorem is lemma 2.2.13
in \cite{Ol1}. We deduce it from the Wick formula.

\proclaim{\smc Theorem (Olshanski)}
Suppose $A_1,\dots,A_k\in\gln\subset\Ugln$. Then
$$
\align
A_1\dots A_k=&
\c A_1\dots A_k\c +\\
&\sum_{1\le a<b\le k} 
\c A_1\dots\pair{A_a\dots A_b}\dots A_k\c +\\
&\sum_{a,b,c}
\c A_1\dots\pair{A_a\dots A_b}
\chainpair{\dots A_c}\dots A_k\c + \dots \,,
\endalign
$$
where the sum is over all (possibly chain) pairings in the set
$\{1,\dots,k\}$.
\endproclaim
For example,
$$
A_1A_2A_3=\c A_1A_2A_3\c + \c\pair{A_1A_2}A_3\c+
\c A_1\pair{A_2A_3}\c+\c\pair{A_1A_2A_3}\c+
\c\pair{A_1A_2}\chainpair{A_3}\c\,.
$$

Note that the sum in theorem is in fact the
sum over all partitions of the set $\{1,\dots,k\}$
into disjoint union of its subsets ({\it clusters\/})
$$
\{i_1,i_2,\dots\},\{j_1,j_2,\dots\},\dots\subset
\{1,\dots,k\}.
$$
Each cluster $\{i_1,i_2,i_3,\dots\}$  corresponds to
the following chain pairing
$$
\c A_1\dots\pair{A_{i_1}\dots A_{i_2}}
\chainpair{\dots A_{i_3}}\lchainpair{\dots\phantom{A_{i_3}}}
\dots A_k\c \,.
$$

\demo{\smc Proof}
Apply the Wick formula to the product
$$
E_{i_1j_1}\dots E_{i_kj_k}
=\sum_{\a_1,\dots,\a_k}
x_{i_1\a_1}\p_{j_1\a_1}
\dots x_{i_k\a_k}\p_{j_k\a_k}\,. \tag 2.4
$$
Remark that the pairing of $\p_{j_p\a_p}$
with $x_{i_q\a_q},p<q$, induces matrix
multiplication of $E_{i_pj_p}$ and $E_{i_qj_q}$. \qed
\enddemo

\head
3. Young basis.
\endhead

\subhead 3.1 \endsubhead
Recall the construction of the Young orthogonal basis in the
irreducible representations of the symmetric groups
$S(k)$, $k=1,2,\dots$. Define it by induction.

The group $S(1)$ is trivial. We can choose any nonzero
vector in its unique irreducible representation. Suppose
$k>1$. Let $\l,|\l|=k$, be a Young diagram and let $V^\l$ be
the corresponding irreducible $S(k)$-module. Let $\mu$
be another Young diagram. Write
$\mu\nea\l$ if $\mu\subset\l$ and $|\mu|=|\l|-1$. The 
Young branching rule asserts that
$$
V^\l=\bigoplus_{\mu\nea\l}V^\mu
\qquad\text{as a $S(k-1)$-module\,.}\tag 3.1
$$
Here the sum is orthogonal with respect to the
$S(k)$-invariant inner product $(\,,\,)$ in $V^\l$.
By definition, the Young basis in $V^\l$ is the union
of the Young bases in direct summands in \tht{3.1}.

It is clear that the Young basis in $V^\l$ is indexed
by the following chains of diagrams
$$
\varnothing=\l^{(0)}\nea\l^{(1)}\nea\dots
\nea\l^{(k-1)}\nea\l^{(k)}=\l\,.
$$
Such a chain is the protocol of a Young
diagram growth from the empty diagram to the diagram $\l$.
This growth can be also represented as follows:
for all $i=1,\dots,k$ put the number $i$ into the
box $\l^{(i)}/\l^{(i-1)}$ of the diagram $\l$. Then
we obtain a {\it Young tableau} of shape $\l$, that is
a tableau $T$ whose entries strictly increase along
each row and down each column. Denote by $v_T$
the Young basis vector corresponding to the tableau $T$.

By our definition each basis vector is defined only
up to a scalar factor. In the sequel we suppose that
$$
(v_T,v_T)=1\,.
$$
This normalization is the only object in this paper
which is not defined over the field ${\Bbb Q}$ of
rational numbers.

Suppose $\a=(i,j)$ is a box of $\l$. Recall that the number
$$
c(\a)=j-i
$$
is called the {\it content} of the box $\a$. For all
$i=1,\dots,k$ put
$$
c_T(i)=c(\l^{(i)}/\l^{(i-1)})\,,
$$
this is the content of the $i$-th box in the tableau $T$.
Observe that always
$$
c_T(1)=0\,.
$$

\subhead 3.2 \endsubhead
For all $i=1,\dots,k$ consider the following
elements of $\RSk$
$$
X_i=(1\,i)+(2\,i)+\dots+(i-1\,\,i)\,. \tag 3.2
$$
In particular, $X_1=0$. These elements were 
introduced by Jucys \cite{Ju} and Murphy \cite{Mu}.
The following proposition is also due to these
authors. Our proof follows \cite{Ol2}, section 4.6.

\proclaim{\smc Proposition} For all $i=1,\dots,k$
$$
X_i v_T=c_T(i) v_T \,.
$$
\endproclaim

\demo{\smc Proof}
For all $p=1,\dots,k$ put
$$
\Si_p =
\sum_{1\le i<j\le p} (i\,j) \quad \in\R[S(p)]\,.
$$
It is clear that $\Si_p$ is a central element
of $\R[S(p)]$ and it is proved, for example, in
\cite{M}, Exercise I.7.7, that in all irreducible
$S(p)$-modules $V^\eta$
$$
\Si_p|_{V^\eta}=
\frac12\sum_i
\left(\eta^2_i-(2i-1)\eta_i\right)\cdot
\id_{V^\eta} \,. \tag 3.3
$$
Clearly,
$$
X_i=\Si_i-\Si_{i-1}\,. \tag 3.4
$$
Choose $q$ so that $\l^{(i)}_q=\l^{(i-1)}_q+1$.
Put $l=\l^{(i)}_q$. By \tht{3.3} and \tht{3.4} we have
$$
\align
X_i|_{V^{\l^{(i-1)}}}&=
\frac12(l^2-(2i-1)l-(l-1)^2+(2i-1)l)
\cdot\id_{V^{\l^{(i-1)}}}\\
&=(l-i)
\cdot\id_{V^{\l^{(i-1)}}}\\
&=c_T(i)
\cdot\id_{V^{\l^{(i-1)}}} \,.
\endalign
$$
Since $v_T\in V^{\l^{(i-1)}}$ this proves the proposition. \qed
\enddemo

Let $T,T'$ be two Young tableaux of shape $\l$. Consider
the matrix element
$$
\pst(s)=(s\cdot v_{T},v_{T'})\,.
$$
Consider the following element of $\RSk$
$$
\Pst=\sum_{s\in S(k)} (s\cdot v_{T},v_{T'})
\cdot s^{-1}\,.
$$

\proclaim{\smc Corollary} For all $i=1,\dots,k$
$$
\align
X_i \Pst &= c_{T}(i) \Pst \tag 3.5\\
\Pst X_i &= c_{T'}(i) \Pst \tag 3.6
\endalign
$$
\endproclaim
\demo{\smc Proof} 
The equalities \tht{3.5} and \tht{3.6} are equivalent to
$$
\align
\sum_{j,j<i} \pst(s(i\,j)) &= c_{T}(i) \pst(s) \tag 3.5'\\
\sum_{j,j<i} \pst((i\,j)s) &= c_{T'}(i) \pst(s)\,, \tag 3.6'
\endalign
$$
which follow from the definition \tht{3.2}, the proposition,
and the invariance of the inner product
$$
(s\cdot v,u)=(v,s^{-1}\cdot u), \quad v,u\in V^\l \,. \qed
$$
\enddemo

It follows from the orthogonality relations for matrix elements
and it also follows from the corollary that 
in the Young basis the operator $\Pst$
is proportional to a matrix unit. Put
$$
P_{TT'} = (\dim\mu\ovr k!) \,\Pst\,,
$$
where $\dim\mu$ is the dimension of $V^\mu$. Then
$$
\alignat2
&P_{TT'}\cdot v_{T'} &&=  v_{T},\\
&P_{TT'}\cdot v_{T''} &&= 0,\quad T''\ne T'\,. 
\endalignat
$$

\remark{\smc Remark}  
The corollary asserts that the matrix elements $\Pst$
form the unique up to scalar factors common eigenbasis
for $2k$ commuting self-adjoint operators which act
by multiplications by $X_1,\dots,X_k$ from the left
and from the right. In fact the representation theory
of the symmetric groups can be rediscovered from some
simple properties of these operators, see \cite{OV}.
\endremark

\subhead 3.3 \endsubhead
The practical computation of matrix elements $\pst$ is a
quite difficult problem. A way of computing them is
the following. First one obtains one particular matrix
element in each irreducible representation and
then the other from this one.

It is known \cite{JK} that in $V^\l$ there is the unique
up to scalar factor vector invariant under the action
of the group
$$
S(\l)=S(\l_1)\times S(\l_2)\times \dots
$$
which is the stabilizer of the subsets
$$
\{1,\dots,\l_1\},
\{\l_1+1,\dots,\l_1+\l_2\},
\{\l_1+\l_2+1,\dots\},\dots \,.
$$
It can be easily deduced from the definition of the Young
basis that this vector is simply the vector $v_T$, where
$T$ is the following Young tableau
$$
T\,=\quad\matrix
1 & 2& \dots& \dots&\dots&\l_1-1&\l_1\\
\l_1+1& \l_1+2&\dots&\l_1+\l_2\\
\l_1+\l_2+1& \dots\\
\dots
\endmatrix\quad,
$$
which is called the {\it row tableau} of shape $\l$. Denote
this tableau by $T^0$. Consider the {\it Young symmetrizer\/}
\cite{JK} corresponding to the tableau $T^0$
$$
\P\Q\,\in\RSk\,,
$$
where
$$
\P=\sum_{s\in S(\l)} s
$$
is the row-symmetrizer of the tableau $T^0$ and
$$
\Q=\sum_{s \text{ preserves columns of } T^0} \sgn(s)\cdot s
$$
is the column-antisymmetrizer of $T^0$.
Both $\P$ and $\Q$ are up to a scalar factor orthogonal
projections
$$
\P^*=\P,\quad \P^2=\l!\,\P,\quad
\Q^*=\Q,\quad \Q^2=(\l')!\,\Q\,,
$$
where $*$ is the involution in $\RSk$ induced by
$$
s\mapsto s^{-1}, \quad s\in S(k)\,,
$$
and $\l!\,=\l_1!\,\l_2!\,\dots\,$. 
The element $\P$ acts in $V^\l$ up to a scalar as the orthogonal
projection onto $v_{T^0}$.

It is known that the product $\P\Q$ vanishes
in any irreducible representation of
$S(k)$ different from $V^\l$. The operator
$$
\P\Q\P=\frac1{(\l')!} \P\Q(\P\Q)^*
$$
is a nonzero operator proportional to the
orthogonal projection onto $v_{T^0}$ and hence
it is proportional to $\Psi_{T^0T^0}$.
It can be easily shown that 
$$
\Psi_{T^0T^0}=\frac1{\l!}\P\Q\P\,.
$$

\subhead 3.4 \endsubhead
Now suppose $v$ is a common eigenvector of the elements $X_1,\dots,X_k$
in a $S(k)$-module $V$
$$
X_i\cdot v = a_i v, \quad i=1,\dots,k,\quad a_i\in\R \,.
$$
There is a standard general method to construct new eigenvectors 
of $X_1,\dots,X_k$ from $v$. Put
$$
s_i= (i\,,i+1), \quad i=1,\dots,k-1 \,.
$$
It is easy to check \cite{Mu} that
$$
\align
s_i X_i +1 &=X_{i+1} s_i\,, \tag 3.7\\
s_i X_j &= X_j s_i\,,\quad j\ne i,i+1\,. \tag 3.8
\endalign
$$
Suppose
$$
a_{p+1}\ne a_p\pm 1 \tag 3.9
$$
for some $p$. Put
$$
v'=\left(
s_p-\frac1{a_{p+1}-a_p} \right)\cdot v \quad \in V\,.
$$
By \tht{3.9} we have $v'\ne 0$. It follows from \tht{3.7} and \tht{3.8}
that
$$
X_i\cdot v'=a_{s_p(i)} v' \,.
$$
It is easy to see that all eigenvectors $\Pst$ in
the $S(k)\times S(k)$-module $\RSk$ can be obtained
in this way from an arbitrary initial matrix element
(for example, $\Psi_{T^0T^0}$) in each irreducible
representation.

\head
4. Proof of the main theorem
\endhead

\subhead 4.1 \endsubhead
We have to prove the matrix equality \tht{1.3}. Prove
that all matrix elements are equal.  
Put $\psi(s)=\pst(s)$ and put $c(i)=c_T(i)$. By \tht{2.2}
we have to prove that for all collections of
indexes 
$$
i_1,\dots,i_k,\quad j_1,\dots,j_k
$$
we have
$$
\multline
\sum_{s\in S(k)} 
\psi(s) \cdot
\left(E_{i_1j_{s(1)}}-c(1)\d_{i_1j_{s(1)}}\right)\dots
\left(E_{i_kj_{s(k)}}-c(k)\d_{i_kj_{s(k)}}\right)
\\ =
\sum_{s\in S(k)} 
\c \psi(s) 
E_{i_1j_{s(1)}}\dots E_{i_kj_{s(k)}}\c\quad.
\endmultline \tag 4.1
$$
To simplify notation put
$$
l_p=j_{s(p)}\,,\quad p=1,\dots,k\,.
$$
The indexes $l_1,\dots,l_k$ vary simultaneously with the
permutation $s\in S(k)$. We are going expand out all
brackets in the LHS of \tht{4.1} and then apply the theorem
from paragraph 2.2 to all monomials in $E_{ij}$.

Fix some $s$ to see what happens. We have the product
$$
\left(E_{i_1l_1}-c(1)\d_{i_1l_1}\right)\dots
\left(E_{i_kl_k}-c(k)\d_{i_kl_k}\right)\,.
$$
First for all $p=1,\dots,k$ we have to choose
in the $p$-th bracket either $E_{i_pl_p}$ or
$(-c(p)\d_{i_pl_p})$. Let us depict our
choice as a diagram like
$$
\ov1\to\circ\qquad
\ov2\to*\qquad
\ov3\to\circ\qquad
\cdots\qquad
\ov{k-1}\to\circ\qquad
\ov{k}\to*\quad,
$$
where the circles represent the factors $E_{i_pl_p}$ and
the asterisks represent the factors $(-c(p)\d_{i_pl_p})$. 
For example, the diagram
$$
\ov1\to\circ\qquad
\ov2\to\circ\qquad
\ov3\to\circ\qquad
\cdots\qquad
\ov{k-1}\to\circ\qquad
\ov{k}\to\circ \tag 4.2
$$
corresponds to the product
$$
E_{i_1l_1}\dots E_{i_kl_k}\,,
$$
and the diagram
$$
\ov1\to*\qquad
\ov2\to*\qquad
\ov3\to*\qquad
\cdots\qquad
\ov{k-1}\to*\qquad
\ov{k}\to* 
$$
corresponds to
$$
\left(-c(1)\d_{i_1l_1}\right)\dots
\left(-c(k)\d_{i_kl_k}\right)\,.
$$
Next, we have to divide the factors $E_{i_pl_p}$
(or, equivalently, the circles in the diagram)
into clusters in all possible ways. This
will be depicted as follows: a cluster
$\{a,b,c\}$
$$
\circ\qquad
\pair{
\ov a\to\circ\qquad
*\qquad
\circ\qquad
\ov b\to\circ}
\chainpair{\qquad
\circ\qquad
\ov c\to\circ}
\qquad *
$$
corresponds to the factor
$$
\c\dots\d_{l_ai_b}\d_{l_bi_c} E_{i_al_c}\dots\c\quad.
$$

We see that the summands which arise in the LHS of \tht{4.1}
are indexed by permutations $s$ and diagrams like
$$
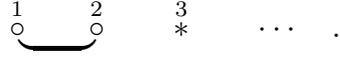
\pair{\ov1\to\circ\qquad\ov2\to\circ}\qquad\overset3\to*
\qquad\cdots\quad.
$$
Denote the corresponding summand by
$$
\bigg[\quad s\quad\bigg|\quad
\pair{\ov1\to\circ\qquad\ov2\to\circ}\qquad\overset3\to*
\qquad\cdots\quad\bigg]\quad.
$$
In order to establish \tht{4.1} we have to show that
all summands cancel each other except those corresponding to
the trivial diagram \tht{4.2}.

\subhead 4.2 \endsubhead
To explain the idea of the proof, we show first that all
summands that contain exactly $k-1$ factors $E_{ij}$
cancel. Such summands correspond to two kind of
diagrams:
$$
\circ\quad\cdots\quad\circ\qquad
\ov b\to\ast\qquad\circ
\quad\cdots\quad\circ\,\qquad b=1,\dots,k\,, \tag 4.3
$$
and
$$
\circ\quad\cdots\quad
\pair{
\ov a\to\circ
\quad\cdots\quad
\ov b\to\circ}
\quad\cdots\quad\circ\,\qquad 1\le a<b\le k\,. \tag 4.4
$$
We claim that for all $s$ and all $b$
$$
\bigg[\quad s\quad\bigg|\,
\cdots\quad
\overset b\to\ast\quad\cdots
\,\bigg]+
\sum_{a<b}\bigg[\,s(ab)\,\bigg|\,
\cdots\quad
\pair{\overset a\to\circ\quad\cdots\quad
\overset b\to\circ}\quad\cdots
\,\bigg]=0\,. \tag 4.5
$$
In fact, all summands in \tht{4.5} are proportional
to
$$
\c \d_{i_bl_b} \prod_{p\ne b} E_{i_pl_p} \c\quad,
$$
and the coefficient equals
$$
-c(b)\psi(s)+\sum_{a,a<b}\psi(s(ab)) \,. \tag 4.6
$$
By \tht{3.5'} this number equals zero. Evidently, by
\tht{4.5} all summands with diagrams \tht{4.3} and \tht{4.4}
cancel.


\subhead 4.3 \endsubhead
Now consider the general case. Suppose we have a
diagram $\G$, for example
$$
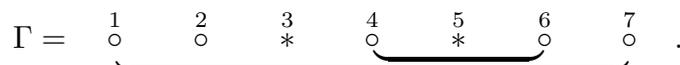
\G=\quad
\dpair{
\overset1\to\circ\qquad
\overset2\to\circ\qquad
\overset3\to\ast\qquad
\pair{
\overset4\to\circ\qquad
\overset5\to\ast\qquad
\overset6\to\circ}\qquad
\overset7\to\circ}\quad.
$$
This diagram corresponds to three clusters
$$
\{1,7\},\{2\},\{4,6\},
$$
and the  subset
$$
\{3,5\}
$$
formed by all asterisks. Denote this asterisk subset by $As$.

Let $b$ be the smallest positive integer such that $b$ is not
a beginning of a new circle cluster. In our example $b=3$. We
say that the diagram $\G$ is of the {\it first kind} if
$b\in As$ and of the {\it second kind} otherwise. For example,
our diagram in example and all diagrams \tht{4.3} are of the
first kind, whereas all diagrams \tht{4.4} are of the second kind.

We claim that for all $s$ and for all diagrams $\G$ of the first
kind the corresponding summand
$$
\bigg[\quad s\quad\bigg|\,
\cdots\quad
\overset b\to\ast\quad\cdots
\,\bigg] \tag 4.7
$$
cancels with the sum
$$
\sum_{a<b}\bigg[\,s(ab)\,\bigg|\,
\cdots\quad
\pair{\overset a\to\circ\quad\cdots\quad
\overset b\to\circ}\quad\cdots
\,\bigg]\,, \tag 4.8
$$
where the pairing of $a$ and $b$ means that $b$ should be
added to the cluster that begins with $a$. 
Note that the diagrams in \tht{4.8} are of the second kind and
all summands with a second kind diagram appear exactly
one time in the sum \tht{4.8} while $s$ ranges over $S(k)$
and $\G$ ranges over all diagrams of the first kind.

In our example
$$
\alignat3
&\bigg[\quad s &&\bigg|
\quad
\dpair{
\overset1\to\circ\qquad
\overset2\to\circ\qquad
\overset3\to\ast\qquad
\pair{
\overset4\to\circ\qquad
\overset5\to\ast\qquad
\overset6\to\circ}\qquad
\overset7\to\circ}\quad
&&\bigg]+\\
&\bigg[\,\, s(1\,3)\,\,&&\bigg|
\quad
\pair{
\overset1\to\circ\qquad
\overset2\to\circ\qquad
\overset3\to\circ}
\dchainpair{\qquad
\pair{
\overset4\to\circ\qquad
\overset5\to\ast\qquad
\overset6\to\circ}\qquad
\overset7\to\circ}\quad
&&\bigg]+\\
&\bigg[\,\,s(2\,3)\,\,&&\bigg|
\quad
\dpair{
\overset1\to\circ\qquad
\pair{\overset2\to\circ\qquad
\overset3\to\circ}
\qquad
\pair{
\overset4\to\circ\qquad
\overset5\to\ast\qquad
\overset6\to\circ}\qquad
\overset7\to\circ}\quad
&&\bigg]=0
\endalignat
$$ 

The cancellation of \tht{4.7} and \tht{4.8} is proved in the
same way as \tht{4.5}. 
It is easy to see that all summands in \tht{4.7} and \tht{4.8}
are proportional. Indeed, suppose 
$\{a,c,d,\dots,z\}$ is the cluster in $\G$ that begins with $a$.
$$
\cdots\,
\pair{
\ov a\to\circ
\,\cdots\,
\ov b\to\ast
\,\cdots\,
\ov c\to\circ}
\chainpair{
\,\cdots\,
\ov d\to\circ}
\lchainpair{
\,\cdots\,}
\,\cdots\,
\rpair{
\,\cdots\,
\ov z\to\circ}
\,\cdots
$$
Then the contribution of this cluster 
to \tht{4.7} is the following factor
$$
\d_{l_ai_c}\d_{l_ci_d}\dots E_{i_al_z} \,.
$$
The contribution of the asterisk on the $b$-th place
is the factor
$$
-c(b) \d_{l_bi_b} \,.
$$
On the other hand, the contribuition of the cluster
$\{a,b,c,d,\dots,z\}$ 
$$
\cdots\,
\pair{
\ov a\to\circ
\,\cdots\,
\ov b\to\circ}
\chainpair{
\,\cdots\,
\ov c\to\circ}
\chainpair{
\,\cdots\,
\ov d\to\circ}
\lchainpair{
\,\cdots\,}
\,\cdots\,
\rpair{
\,\cdots\,
\ov z\to\circ}
\,\cdots
$$
to the $a$-th summand in \tht{4.8} is the factor
$$
\d_{l_bi_b}\d_{l_ai_c}\d_{l_ci_d}\dots E_{i_al_z} \,.
$$
Therefore all summands in \tht{4.7} and \tht{4.8} are
proportional. The coefficient equals \tht{4.6} again and hence
equals zero. This concludes the proof of the theorem. 


\head
5. Quantum immanants and higher Capelli identities.
\endhead

\subhead 5.1 \endsubhead
In this section we specialize the main theorem 
for central elements of $\Ugln$. 
Recall that the trace of an element of 
$\Ugln\otimes \M(n)^{\otimes n}$ is defined by
$$
\tr \big( \sum_
{i_1,j_1,\dots,i_n,j_n}
A_{i_1,j_1,\dots,i_n,j_n} \otimes
e_{i_1,j_1}\otimes\dots\otimes e_{i_n,j_n} \big) =
\sum_{i_1\dots,i_n}
A_{i_1,i_1,\dots,i_n,i_n} \, \in \Ugln \,.
$$
Let $\mu,|\mu|=k$ be a Young diagram. Denote by
$\Tab(\mu)$ the set of all Young tableaux of shape $\mu$.
Put
$$
\dim\mu=\dim V^\mu\,.
$$
Recall that we consider the character 
$\chm$ of the module $V^\mu$
as an element of $\RSk$
$$
\chm = \sum_{s\in S(k)} \chm(s) \cdot s \, \in \RSk \,.
$$ 
Put
$$
\ET={\Bbb E}_{TT}\,,
$$
where the fusion matrix $\ETT$ is defined by 
$$
\ETT=(E-c_T(1))\otimes \dots \otimes (E-c_T(k))
\cdot \Pst \,.
$$
The following theorem is the main result of \cite{Ok}.
Here we deduce its first claim from the main theorem.

\proclaim{\smc Theorem \cite{Ok}}
\roster
\item"(a)"
For all $T\in\Tab(\mu)$
$$
\tr \ET =
\frac1{\dim\mu}\, \tr X^{\otimes k}
(D')^{\otimes k} \chm \quad\in \Ugln \,. \tag 5.1
$$
In particular, the LHS of \tht{5.1} does not depend
on the choice of $T\in\Tab(\mu)$.
\item"(b)"
The element \tht{5.1} lies in the center $\Zgln$ of
$\Ugln$.
\item"(c)"
The elements \tht{5.1} form a linear basis of $\Zgln$
indexed by all Young diagrams $\mu$.
\endroster
\endproclaim
\demo{\smc Proof}
Prove (a). By the main theorem we have
$$
\tr \ET =
\tr X^{\otimes k}
(D')^{\otimes k} \Psi_T \,. 
$$
Since the entries of the matrix $X$ commute we have
$$
s\cdot X^{\otimes k} =
 X^{\otimes k}\cdot s\,, \tag 5.2
$$
for all $s\in S(k)$, and similarly
$$
s\cdot (D')^{\otimes k} =
(D')^{\otimes k} \cdot s\,. \tag 5.3
$$
Observe that
$$
1\ovr k! \sum_{s\in S(k)} s\Psi_T s^{-1} =
\frac1{\dim\mu}\, \chm \,.
$$
Therefore
$$
\align
\tr \ET &=
1\ovr k!\sum_{s\in S(k)}\tr s\, X^{\otimes k}
(D')^{\otimes k} \Psi_T \,s^{-1} \\
&= 
\frac1{\dim\mu}\, \tr X^{\otimes k}
(D')^{\otimes k} \chm \,.
\endalign
$$

Prove (b). Prove, for example, that the LHS of \tht{5.1}
is a central element.
Denote by $g_{ij}$ and $\gi_{ij}$ the matrix elements
of a matrix $g\in GL(n)$ and its inverse matrix $g^{-1}$.
The following equality is obvious
$$
\sum_k g_{ki} \gi_{jk} = \d_{ij}\,. \tag 5.4
$$
Consider the adjoint action $\Ad(g)$ of $g$ in $\gln$
$$
\Ad(g)\cdot E_{ij} = \sum_{k,l} g_{ki} \gi_{jl} E_{kl} \,. \tag 5.5
$$
Under the adjoint action of $g$ the 
entries of the matrix $(E-u)$ are
transformed as follows
$$
\alignat2
(E-u) @>\Ad(g)>> 
& \sum_{i,j} 
\left(
\sum_{k,l} g_{ki} \gi_{jl} E_{kl} \right) \otimes
e_{ij} - u\sum_{i}1\otimes e_{ii} 
\quad&&\text{by \tht{5.5}} \\
&= \sum_{k,l} 
\left(
E_{kl} - u\d_{kl} \right) \otimes
\left(\sum_{i,j} g_{ki} \gi_{jl} e_{ij} \right) 
\quad&&\text{by \tht{5.4}} \\
&= g' (E-u) (g')^{-1} \tag 5.6
\endalignat
$$
The product \tht{5.5} is the product of the matrix $(E-u)$
with entries in $\Ugln$ and two matrices with entries
in the ground field.
Consider the following element of $\Ugln$
$$
\tr((E-u_1)\otimes\dots\otimes (E-u_k)\cdot s)\,,  \tag 5.7
$$
where the numbers $u_i$ and 
the permutation $s\in S(k)$ are arbitrary. By \tht{5.4} the adjoint action
of $g'$ takes this element of $\Ugln$ to
$$
\tr(g^{\otimes k}
(E-u_1)\otimes\dots\otimes (E-u_k)
(g^{-1})^{\otimes k}
\cdot s) =
\tr((E-u_1)\otimes\dots\otimes (E-u_k)\cdot s)\,.
$$
Hence,\tht{5.7} is an element of $\Zgln$.
Therefore \tht{5.1} is an element of $\Zgln$.

Prove (c). Consider the standard filtration in $\Ugln$
and consider the isomorphism
$$
\gr \Ugln \cong S(\gln)\,.
$$
It is clear that
$$
\tr \ET =
\frac1{\dim\mu}\, \tr E^{\otimes k} \chm +
\text{  lower terms }\,.
$$
Suppose $G=(g_{ij})$ is a $n\times n$-matrix. It
follows from the classical decomposition of the
vector space of tensors that the following polynomial
in $g_{ij}$
$$
\tr\, G^{\otimes k} \chm\ovr k! \tag 5.8
$$
equals the trace of $G$ in the irreducible $GL(n)$-module
with highest weight $\mu$ (or, equivalently, it equals
the Schur polynomial $s_\mu$ in the eigenvalues of $G$).
The polynomials \tht{5.8} form a linear basis in the
vector space of invariants for the adjoint action of $GL(n)$
on $\gln$. Hence the elements \tht{5.1} form a linear
basis in $\Zgln$. \qed
\enddemo

\remark{\smc Remark} Given a matrix $A=(a_{ij})$, $i,j=1,\dots,k$,
the number
$$
\sum_{s\in S(k)} \chm(s) \, a_{1,s(1)} \dots a_{k,s(k)}
$$
is called the {\it $\mu$-immanant} of the matrix A.
If $\mu=(1^k),(k)$ then the $\mu$-immanant turns into
determinant and permanent respectively.
Observe that \tht{5.8} is
the sum of $\mu$-immanants of principal
$k$-submatrices (with repeated rows and columns) of the matrix $G$.
\endremark

\subhead 5.2 \endsubhead
By definition, put
$$
\Sm=\frac{\dim\mu}{k!} \,\tr \ET \,,
\quad T\in\Tab(\mu)\,. \tag 5.9
$$
By the theorem this central element does not depend on the
choice of $T\in\Tab(\mu)$. If $\mu=(1^k)$ then the
definition of $\Sm$ turns into the definition of 
quantum determinant for $\Ugln$ (see \cite{KuS} or \cite{MNO}).
By analogy to quantum determinant and because of the
structure of the highest term of \tht{5.1} we call
$\Sm$ the {\it quantum $\mu$-immanant}. Quantum
immanants were introduced and studied in the 
authors paper \cite{Ok}; from a different
point of view they were studied in \cite{OO}.
Here we mention some most important properties of these remarkable
basis elements of $\Zgln$.

\subhead 5.3 \endsubhead
We claim that the identity \tht{5.1} is a direct
generalization of the classical Capelli identity \tht{1.2}.
If $\mu=(1^k)$ then it is easy to see that the
RHS of \tht{5.1} turns into the RHS of \tht{1.2}. 
Let us show that the LHS of \tht{5.1} turns into
the LHS of \tht{1.2}. Let
$$
i_1,i_2,\dots, i_k
$$ 
be a $k$-tuple of indexes. Denote by $\iota!$ the
order of the stabilizer of this collection in
the symmetric group $S(k)$. 
For example, if all $i_j$ are distinct then $\iota!=1$.
\proclaim{\smc Teorem \cite{Ok}} 
$$
\align
\Sm&=\sum_{i_1\ge\dots\ge i_k} 1/{\iota!} 
\sum_{T\in\Tab(\mu)} \sum_{s\in S(k)}
\psi_T(s) \, (E_{i_1i_{s(1)}})
(E_{i_2i_{s(2)}}-c_T(2)\d_{i_2i_{s(2)}}) \dots \\
&=\sum_{i_1\le\dots\le i_k} 1/{\iota!} 
\sum_{T\in\Tab(\mu)} \sum_{s\in S(k)}
\psi_T(s) \, (E_{i_1i_{s(1)}})
(E_{i_2i_{s(2)}}-c_T(2)\d_{i_2i_{s(2)}}) \dots  
\endalign 
$$
\endproclaim

In \cite{Ok} this theorem was used in proof of the identity \tht{5.1}.
Here we deduce this theorem from \tht{5.1}.

\demo{\smc Proof}
By \tht{5.1} we have
$$
\align
\Sm&=1\ovr k! \sum_{T\in\Tab(\mu)} \tr \ET \\ 
&=1\ovr k! \sum_{i_1,\dots,i_k}  
\sum_{T\in\Tab(\mu)} \sum_{s\in S(k)}
\psi_T(s) \, (E_{i_1i_{s(1)}})
(E_{i_2i_{s(2)}}-c_T(2)\d_{i_2i_{s(2)}}) \dots  \tag 5.10
\endalign 
$$
By \tht{5.2} and \tht{5.3} the  matrix
$$
\sum_{T\in\Tab(\mu)} \ET  =
\tr X^{\otimes k} (D')^{\otimes k} \chm 
$$
is invariant under conjugation by element of the
group $S(k)$. Hence all $k!\ovr\iota!$ different
rearrangement of the indexes $i_1,\dots,i_k$ make the same
contribution to the trace $\tht{5.10}$ and hence
we can choose an arbitrary (for example,
increasing or decreasing) ordering of the indexes.
This proves the theorem. \qed
\enddemo

It is easy to see that the second formula for $\Sm$
in the theorem turns into the LHS of \tht{1.2} when
$\mu=(1^k)$. Therefore we call the equalities \tht{5.1}
the {\it higher Capelli identities}.

\remark{\smc Remark}
The aguments based on the two trivial observations \tht{5.2}
and \tht{5.3} and on the main theorem provide an elementary proof
of many identities involving the  matrix $\ET $
which seemed to require deep machinery
of Yangians and R-matrixes. One of them
is the following identity
$$
\align
P_{T^1T^2}{\Bbb E}_{T^3T^4}&=
P_{T^1T^2} X^{\otimes k} (D')^{\otimes k} \Psi_{T^3T^4}\\
&=X^{\otimes k} (D')^{\otimes k} P_{T^1T^2} \Psi_{T^3T^4}\\
&=X^{\otimes k} (D')^{\otimes k} \d_{T^2T^3} \Psi_{T^1T^4} \\
&= \d_{T^2T^3}\, {\Bbb E}_{T^1T^4} \,,
\endalign
$$
where $T_1,\dots,T_4$ are four arbitrary Young tableaux.
\endremark

\subhead 5.4 \endsubhead
Denote by $\pi_\l$ the representation of $GL(n)$
with highest weight $\l$. Recall the definition
of the {\it shifted Schur function} from \cite{OO}.
Put
$$
(x\f k)= x(x-1)\dots(x-k+1),.
$$
This product is called falling {\it factorial
power}. Put also 
$$
\rho=(n-1,\dots,1,0)\,
$$
By definition 
$$
\sm(x_1,\dots,x_n) =\frac
{\det \big[(x_i+\rho_i\f \mu_j+\rho_j)\big]}
{\det \big[(x_i+\rho_i\f \rho_j)\big]} \,. 
$$
The relation between quantum immanants and
shifted Schur functions is the following:

\proclaim{\smc Theorem \cite{Ok}}
Put $\sm(\l)=\sm(\l_1,\l_2,\dots)$. Then
$$
\pi_\l(\Sm)=\sm(\l)\,. \tag 5.11
$$
\endproclaim

Shifted Schur function have many
remarkable properties \cite{OO}
(see also \cite{Ok}). Most of these
properties have a natural interpretation in
terms of quantum immanants $\Sm$ and
are closely related to higher Capelli
identities.  One of the main technical tools to
handle shifted Schur functions is the
following theorem (we shall need
it in the proof of \tht{5.11}). An argument
very close to our proof was used by S.~Sahi
in \cite{S}. Shifted Schur function are a
particular case of certain remarkable polynomials,
which existence was proved in \cite{S}. This
particular case is much more simple and can
be studied much deeper than the general case 
considered in \cite{S}.

By $\Ls(n)$ denote the algebra of polynomials in
variables $x_1,\dots,x_n$ which are symmetric in
new variables $x_1+\rho_1,\dots,x_n+\rho_n$. 
Such  polynomials are called {\it shifted
symmetric} \cite{OO}. It is clear that $\sm\in\Ls(n)$.
Denote by $H(\mu)$ the product of the hook lenghts of
all boxes of $\mu$.

\proclaim{\smc Characterization theorem \cite{Ok}}
Any of the two following properties determines the
polynomial $\sm\in\Ls(n)$ uniquely:
\roster
\item"(A)"
$\deg {\sm}\le|\mu|$ and
$$
\sm(\l)=\d_{\mu\l} H(\mu)
$$
for all $\l$  such that $|\l|\le|\mu|$;
\item"(B)"
the highest term of $\sm$ is the 
ordinary Schur function $s_\mu$  and
$$
\sm(\l)=0 
$$
for all $\l$ such that $|\l|<|\mu|$.
\endroster
\endproclaim

\demo{\smc Proof of \tht{5.11}}
It is well known that the eighenvalue in the
representation $\pi_\l$ of any element of $\Zgln$
is a shifted symmetric function in $\l$. 
Apply $\Sm$ to the highest vector $v$ in the representation
$\pi_\l$. We have
$$
\alignat2
E_{ii}\cdot v & = \l_i v\,, \qquad && i=1,\dots,n\\
E_{ij}\cdot v & = 0\,,\qquad && i<j \,.
\endalignat
$$
By arguments used in proof of part (c) of the theorem in section 5.1
$$
\pi_\l(\Sm)=s_\mu(\l)+\text{  lower terms  }.
$$
On the other hand it is clear that $\Sm$ vanishes
in all representations $\pi_\l$ such that $|\l|<|\mu|$.
Indeed, these reprentations arise as subrepresentations  
of the representation of $\Ugln$ in the  vector space of 
polynomials on $\M(n,m)$
of degree $|\l|$. Such polynomials are
clearly annihilated by the differential
operator in the RHS of \tht{5.1}. Now \tht{5.11}
follows from the characterization theorem. \qed
\enddemo


\enddocument

\end